\documentclass[10pt]{article}
\usepackage[latin1]{inputenc}
\usepackage{amsmath}
\usepackage{amsfonts}
\usepackage{amssymb}
\usepackage[american]{babel}
\usepackage{graphicx}
\usepackage{amsthm}
\usepackage{pstricks}

\begin{document}

\title{Covering cubic graphs with matchings of large size}

\author{S.Bonvicini \thanks{Dipartimento di Scienze e Metodi dell'Ingegneria,
Universit\`a di Modena e Reggio Emilia, via Amendola 2, 42100 Reggio
Emilia (Italy)}
G.Mazzuoccolo \thanks{G-Scop Laboratory, Grenoble, France. Research supported by the project ``INdAM fellowships in mathematics and/or applications for experienced researchers cofunded by Marie Curie actions''}}

\maketitle

\newtheorem{proposition}{Proposition}
\newtheorem{theorem}{Theorem}
\newtheorem{lemma}{Lemma}
\newtheorem{definition}{Definition}
\newtheorem{example}{Example}
\newtheorem{corollary}{Corollary}
\newtheorem{conjecture}{Conjecture}

\begin{abstract} \noindent
Let $m$ be a positive integer and let $G$ be a cubic graph of order $2n$. 
We consider the problem of covering the edge-set of $G$ with the minimum number of matchings of size $m$. This number is called  excessive $[m]$-index of $G$ in literature. The case $m=n$, that is a covering with perfect matchings, is known to be strictly related to an outstanding conjecture of Berge and Fulkerson. In this
paper we study in some details the case $m=n-1$. We show how this parameter can be large for cubic graphs with low connectivity and we furnish some evidence that each cyclically $4$-connected cubic graph of order $2n$ has excessive $[n-1]$-index at most $4$. Finally, we discuss the relation between excessive $[n-1]$-index and some other graph parameters as oddness and circumference.
\end{abstract}

\noindent \textit{ Keywords: excessive index, Berge-Fulkerson conjecture, matchings, cubic graphs. 
MSC(2010): 05C15, 05C70}

\section{Introduction}\label{sec:intro}

Throughout this paper, a graph $G$ always means a cubic simple
connected finite graph (without loops and parallel edges).
We refer to any introductory book for graph-theoretical notation and terminology not described in this paper (see for instance \cite{BonMur}) .

The excessive $m$-index of a graph $G$, denoted by $\chi'_{[m]}(G)$,
is first defined by Cariolaro and Fu in \cite{CarFu2} as the minimum
number of matchings of size $m$ needed to cover the edge-set of $G$.
In what follows, $[m]$-matching will stand for a matching of size $m$.
The excessive $m$-index of particular classes of graphs is computed
(for some values of $m$) in \cite{CarFu2} and \cite{CarFu3} and a general formula for
small values of $m$ is furnished in \cite{CarFu}. It is recently
proved by the second author (see \cite{Maz} and \cite{Maz2}) that
the excessive $m$-index of a graph is strictly related to the
well-known Berge-Fulkerson conjecture and its generalization given
by Seymour in \cite{Se}. Mainly for this reason we will focus our
attention on cubic graphs. We would like to stress that for small
values of $m$ the problem is already solved, whereas it remains
 open for large values of $m$. More precisely, by a result of
Cariolaro, if $\left\lceil \frac{|E(G)|}{m} \right\rceil > \chi'(G)$ holds,
then $\chi'_{[m]}(G)=\left\lceil \frac{|E(G)|}{m} \right\rceil$. Hence,
by direct computation the following proposition holds:
\begin{proposition}\label{smallmatchings}
Let $G$ be a cubic graph of order $2n$. If $m<\left\lceil \frac{3n}{4}
\right\rceil$, then $\chi'_{[m]}(G)=\left\lceil \frac{3n}{m} \right\rceil$.
\end{proposition}
That naturally leads our attention to large values of $m$.
The largest possible case, that is $m=n$, is completely open: it is
conjectured (see Conjecture \ref{BergeFulkerson}) that $\chi'_{[n]}(G)
\leq 5$ for each $2$-connected cubic graph $G$, but it is still
unproved the existence of a constant $k$ such that $\chi'_n(G) \leq
k$ for each $2$-connected cubic graph $G$.
To the best of our knowledge, the best upper bound for the excessive
$n$-index of a $2$-connected cubic graph in terms of its size is given in \cite{Maz3}.\\

In the present paper we will focus our attention on the case $m=n-1$. 
In particular, we will address the following question:\\

\textit{ How many $[n-1]$-matchings we need to cover
the edge-set of a cubic graph of order $2n$?} \\

\begin{figure}[h]
\centering
\includegraphics[scale=0.2]{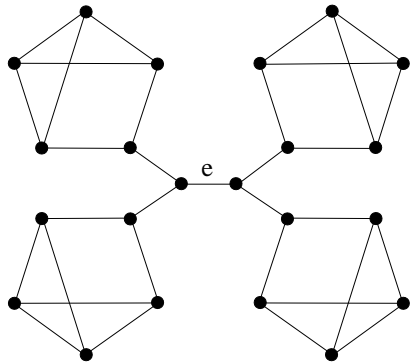}
\caption{A graph with no $[n-1]$-covering}\label{infinito}
\end{figure}

Trivially, if a cubic graph $G$ has an edge which is contained in no
$[n-1]$--matching, then it is not possible to cover the edge--set of
$G$ with $[n-1]$--matchings, in this case we set
$\chi'_{[n-1]}(G)=\infty$. For instance, the graph in Figure
\ref{infinito} is a graph of order $22$ with $\chi'_{[10]}(G)=\infty$:
one can easily check that the edge labelled $e$ does not belong to a
$[10]$-matching of $G$.

In Section $2$, we show that even under the assumption that the excessive $[n-1]$-index of $G$ is finite there exists a family of
cubic graphs with arbitrary large excessive $[n-1]$-index. In Section $3$, we prove some general lemmas and propositions that we will use in the proofs of last section.
In Section $4$, we consider $2$-connected cubic graphs. First of all, by assuming the Berge-Fulkerson conjecture true, we obtain that $\chi'_{[n-1]}(G)$ is at most five for each $2$-connected cubic graph $G$. In particular, we provide examples of $2$-connected and $3$-connected cubic graphs with excessive $[n-1]$-index equals to five. After that, we prove that $\chi'_{[n-1]}(G)\leq 4$ holds for cubic graphs of oddness at most $4$, for cubic graphs with circumference at least $2n-2$ and for the class of $3^*$-connected cubic graphs introduced by Albert, Aldred, Holton and Sheehan in \cite{AlbAldHolShe}.

\section{$1$-connected cubic graphs}

In this section, we consider $1$-connected cubic graphs such
that every edge of the graph is contained in at least one
$[n-1]$-matching (where $2n$ is the order of the graph).

In the next proposition, we show that there exist graphs with finite
excessive $[n-1]$-index as large as we want.

\begin{proposition}
Let $m\geq 3$, there exists a cubic graph $G$ of order $2n$ such
that $\chi'_{[n-1]}(G)\geq m$.
\end{proposition}

\begin{figure}[h]
\centering
\includegraphics[width=5cm]{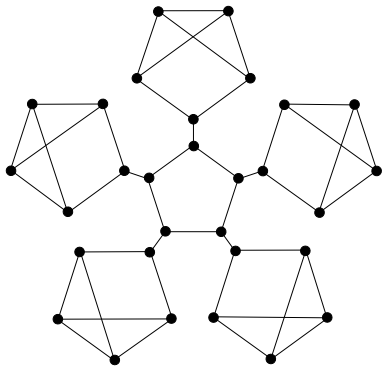}
\caption{A graph of order $30$ with $\chi'_{[14]}=5$.} \label{1connesso}
\end{figure}

\textit{Proof.} Let $C_m=(v_1,\ldots, v_m)$ be a cycle of length
$m\geq 3$. For $i=1,\ldots, m$, denote by $G_i$ a connected graph
sharing no vertex with $C_m$, having $2m_i$ vertices of degree $3$
and exactly one vertex of degree $2$, say $u_i$. Furthermore, $G_i$ is such that each
edge is contained in at least one $[m_i]$-matching of $G_i$. Let
$G$ be the graph having $V(G)=V(C_m)\cup V(G_1)\cup\ldots\cup
V(G_m)$ as vertex--set and $E(G)=E(C_m)\cup E(G_1)\cup\ldots\cup
E(G_m)\cup\{\{u_i, v_i\} | 1\leq i\leq m\}$ as edge--set. The edges
$\{u_i, v_i\}$'s are bridges and will be called \emph{spokes} (in
Figure \ref{1connesso}, we depicted an example with $30$ vertices).
We set $2n=|V(G)|=2m+\sum^m_{i=1} m_i$ and show that no
$[n-1]$-matching of $G$ can contain two edges of $C_m$.

Assume there exists an  $[n-1]$-matching of $G$, say $M$,
containing $h\geq 2$ edges of $C_m$, then $M$ can contain at most
$m-2h$ spokes. Since $M$ can contain at most $m_i/2$ edges of $G_i$,
we have $|M|\leq\sum^m_{i=1} m_i/2+(m-2h)+h=n-h\leq n-2$, a
contradiction, since $M$ has size $n-1$. Hence, each
$[n-1]$-matching of $G$ contains at most one edge of $C_m$. Since
$|E(C_m)|=m$, we need at least $m$ $[n-1]$-matching to cover the
edges of $G$.\qed

\section{General properties}

\begin{lemma}\label{lemmacovering}
Let $G$ be a cubic graph of order $2n$.
Let ${\cal M}=\{M_1,\ldots,M_t\}$ be a covering of $G$ such that
$\sum_{i=1}^t |M_i|=(n-1)t$. Then, there exists a $[n-1]$-covering
of $G$ of size $t$.
\end{lemma}

\textit{Proof.} Suppose $\cal M$ is not a $[n-1]$-covering of $G$,
then there exist two matchings $M_i$ and $M_j$ such that
$|M_i|<n-1<|M_j|$. Consider the subgraph $M_i \cup M_j$, $i\neq j$,
of $G$: a connected component of $M_i \cup M_j$ is either a path or
a cycle of even length (eventually a single edge belonging both to
$M_i$ and $M_j$). In the latter case the connected component of $M_i
\cup M_j$ has the same number of edges of $M_i$ and $M_j$. Then, from
$|M_i|<n-1<|M_j|$, there exists at least a connected component
consisting of a path $P$ of odd length (in case a single edge) starting and finishing with
edges of $M_j$. The exchange of edges in $P$ increases once $|M_i|$ and decreases once $|M_j|$. The iteration
of this process furnishes a $[n-1]$-covering of $G$ of size $t$. $\qed$ \\

\begin{proposition}\label{class1}
Let $G$ be a $3$--edge--colorable cubic graph of order $2n\geq 8$.
Then $\chi'_{[n-1]}(G)=4$.
\end{proposition}

\textit{Proof.} Let $M_1$, $M_2$ and $M_3$ be three pairwise
disjoint perfect matchings of $G$. Let $M_4\subset M_1$ be an arbitrary
$[n-4]$-matching of $G$. The set $\{M_1,\ldots, M_4\}$ is a
covering of $G$ such that $\sum_{i=1}^4 |M_i|=4(n-1)$.
The assertion follows from Lemma \ref{lemmacovering}.  $\qed$ \\

\smallskip

There are only three cubic graphs of order $2n<8$ which are
$3$--edge--colorable, namely, the complete graph $K_4$, the complete
bipartite graph $K_{3,3}$ and the prism $Y_3$ on $6$ vertices. By Proposition \ref{smallmatchings}, $\chi'_{[n-1]}(K_4)=6$ (as $n-1=1$ and
$K_4$ has $6$ edges), $\chi'_{[n-1]}(K_{3,3})=\chi'_{[n-1]}(Y_3)=5$
(as $n-1=3$, $K_{3,3}$ and $Y_3$ have $9$ edges).

\begin{lemma}\label{MN}
Let $G$ be a cubic graph of order $2n$, with $2n \geq 8$. If there
exist a perfect matching $M$ and a $[n-1]$-matching $N$ of $G$ with
empty intersection, then $\chi'_{[n-1]}(G)=4$.
\end{lemma}
\textit{Proof.} Denote by $H$ the complementary subgraph of $M \cup
N$ in $G$. The subgraph $H$ has all vertices of degree one but two
vertices, say $u$ and $v$, of degree two. If the vertices $u$ and
$v$ are adjacent in $G$, then $N \cup \{[u,v]\}$ is a perfect
matching of $G$ disjoint from $M$. Since a cubic graph with two
disjoint perfect matchings is $3$-edge-colorable, the assertion
follows from Proposition \ref{class1}. Now consider the case $u$ and
$v$ non-adjacent in $G$: the set of edges with both endvertices of
degree $1$ in $H$ is a $[n-3]$-matching, say $L$, of $G$. Denote by
$e_1$, $e_2$ (respectively $f_1$,$f_2$) the edges of $H$ incident $u$
(respectively $v$). Set $L_1=L \cup \{ e_1, f_1\}$ and $L_2=(L \cup
\{e_2,f_2\}) \setminus \{ e \}$, where $e$ is an arbitrary edge of
$L$ (such an edge does exist by the assumption on the order of $G$).
The set $\{M,N,L_1, L_2\}$ satisfies Lemma \ref{lemmacovering}, hence the assertion follows. $\qed$ \\

\section{$3$-graphs}

A cubic graph is a $3$--graph if and only if it is bridgeless. We
recall that an $r$--graph is an $r$--regular graph $G$ of even order
such that every edge--cut which separates $V(G)$ into two sets of
odd cardinality has size at least $r$. This notion was introduced in
\cite{Se}. An $r$--graph $G$ is $1$--extendable (every edge of $G$
is contained in a perfect matching), hence
$\chi'_{[n-1]}(G)<\infty$.

It is recently proved by the second author (see \cite{Maz}) that the
well-known conjecture of Berge and Fulkerson \cite{Ful} can be stated as
follows:

\begin{conjecture}[Berge-Fulkerson]\label{BergeFulkerson}
Let $G$ be a $3$--graph. Then, $\chi'_{[n]}(G) \leq 5$.
\end{conjecture}

In the next proposition, we show how the Berge-Fulkerson conjecture
implies that also the excessive $[n-1]$-index of a $3$--graph is
bounded by a constant.

\begin{proposition}\label{BergeFulkersonimplication}
The Berge-Fulkerson conjecture implies $\chi'_{[n-1]}(G)\leq 5$ for
each $3$-graph $G$ of order greater than $4$.
\end{proposition}

\textit{Proof.} Let $G$ be a $3$-graph. If $G$ is
$3$--edge--colorable, then the result follows from Proposition
\ref{class1}. Consider $G$ non $3$-edge-colorable. By assuming that the
Berge-Fulkerson conjecture is true, we have $5$ perfect matchings of
$G$ covering the edge--set of $G$ such that each edge belongs to at
most two perfect matchings. The intersection of each pair of these
perfect matchings is non empty, since the existence of two disjoint perfect matchings implies the existence of a $3$-edge-coloring. Hence, intersecting in pairs the $5$ perfect matchings, we find at least
$10$ distinct edges of $G$ belonging to exactly two of the five
perfect matchings. We can delete five of these edges once from the
perfect matchings. We obtain a covering of cardinality $5n-5$ and
the assertion follows from
Lemma \ref{lemmacovering}. $\qed$ \\

 By using an analogous argument, it is possible to prove that cubic graphs with excessive $[n]$-index $4$ have excessive $[n-1]$-index equals to $4$ too. It is proved in \cite{FouVan2} that some classical families of snarks (i.e. non $3$-edge-colorable cubic graph with girth at least $5$ and cyclically $4$-connected), as Flower snarks and Blanusa snarks,  
have excessive $[n]$-index $4$. Hence, we can state that the excessive $[n-1]$-index is $4$ for snarks of these families.\\
A generation of all snarks up to $36$ vertices is performed by Brinkmann, Goedgebeur, H\"agglund and Markstr\"om (see \cite{Bri}). Janos H\"agglund verified that all snarks with at most $32$ vertices have two perfect matchings $M_1$ and $M_2$ such that their intersection is a unique edge $e$. Since $M_1$ and $M_2 \setminus \{ e \}$ satisfies the hypothesis of Lemma \ref{MN} we can state the following proposition:

\begin{proposition}
Let $G$ be a snark of order $2n \leq 32$. Then, $\chi_{[n-1]}(G)=4$.
\end{proposition}

In what follows we will construct examples of $2$-connected and
$3$-connected cubic graphs having $\chi'_{[n-1]}(G)>4$, furthermore
we will give some evidence to the fact that the excessive
$[n-1]$-index of a cyclically $4$-connected cubic graph could be at
most $4$.

\subsection{$2$-connected}

\begin{figure}[h]
\centering
\includegraphics[width=6cm]{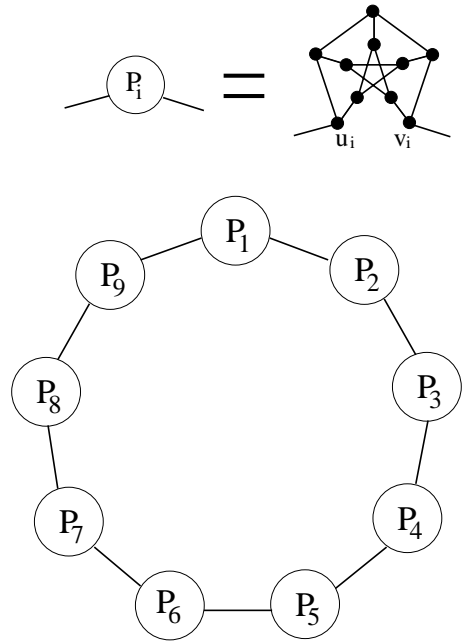}
\caption{A $2$-connected cubic graph $G$ of order $90$ with $\chi'_{[44]}(G)>4$} \label{2connesso}
\end{figure}

In this section, we furnish an example of a $3$-graph $G$ of order $2n$ with $\chi'_{[n-1]}(G)>4$.
Let $G$ be the graph obtained in the following way: consider nine copies $P_i$, for $i=0,\ldots,8$, of the graph obtained by removing from the Petersen graph an edge $[u_i,v_i]$. Add the edges $[v_i,u_{i+1}]$, indices taken modulo $9$ (see Figure \ref{2connesso}).

Suppose there exists a $[n-1]$--covering $\mathcal M=\{M_1,M_2,M_3,M_4\}$ of $G$ of size four. Each $M_i$ leaves two vertices of $G$ uncovered, so we have at most eight vertices which are uncovered at least once. It follows that in at least one of the nine $P_i$'s, without loss of generality we can suppose $P_0$, all vertices are covered by each $[n-1]$-matching. Starting from $P_0$, we  obtain a copy of the Petersen graph by removing the two edges $[u_0,v_8]$ and $[v_0,u_1]$ and adding the edge $[u_0,v_0]$. \\
Consider each $[n-1]$-matching $M_j$: both the edges $[u_0,v_8]$ and $[u_0,v_1]$ belong to $M_j$ or neither one does. In the first case, we add to $M_j \cap P_0$ the edge $[u_0,v_0]$ in order to obtain a perfect matching of $P$, in the latter case $M_j \cap P_0$ is a perfect matching of $P$. This implies that we have a covering of the Petersen graph with four perfect matchings. Since it is well-known that $\chi'_{[n]}(P)=5$ (see for istance \cite{BonCar}), we have a contradiction. \\
Then, the excessive $[n-1]$-index of $G$ is at least $5$ (one can check that it is actually $5$). 

\subsection{$3$-connected}

In this section, we construct an example of a $3$-connected cubic
graph $G$ of order $2n$ with $\chi'_{[n-1]}(G)>4$.

\begin{figure}[h]
\centering
\includegraphics[width=5cm]{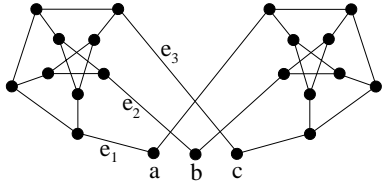}
\caption{The graph $H$} \label{H}
\end{figure}

\begin{figure}[h]
\centering
\includegraphics[width=7.5cm]{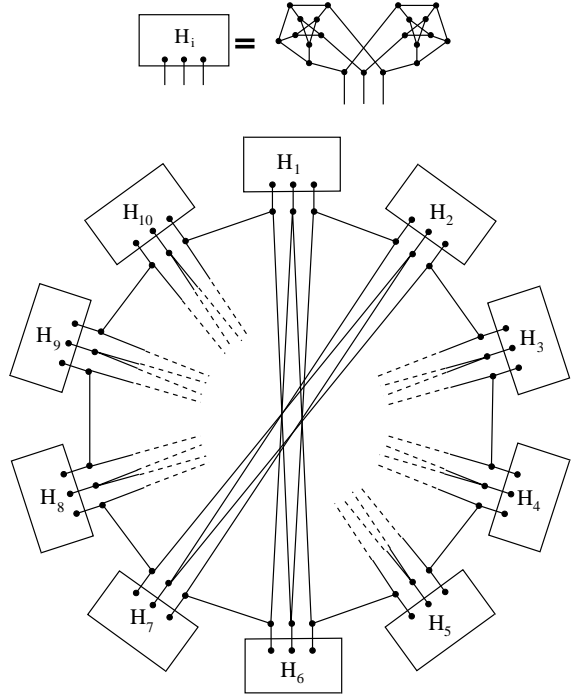}
\caption{A $3$--edge connected cubic graph $G$ of order
$240$, with $\chi'_{119}(G)>4$.} \label{3connesso}
\end{figure}

Consider the graph $H$ in Figure \ref{H}.
The graph $H$ is obtained starting from two copies of the Petersen graph $P$: in
each copy of $P$ delete a vertex, add three
new vertices, say $a$, $b$, $c$, and construct the edges $e_i$,
$e'_i$, $i=1$, $2$, $3$. To construct the graph $G$, consider $2m$
copies of $H$, say $H_1,\ldots, H_{2m}$; for $1\leq i\leq 2m$, label
the vertices of degree $2$ in $H_i$ by $a_i$, $b_i$, $c_i$ and add
three new vertices $u_i$, $v_i$, $w_i$. Construct the edges $[a_i,
u_i]$, $[b_i, v_i]$, $[c_i, w_i]$, with $1\leq i\leq 2m$; $[v_i,
u_{i+m}]$, $[v_i, w_{i+m}]$, with $1\leq i\leq m$; $[u_i, w_{i-1}]$,
with $1\leq i\leq 2m$ (the subscripts are read modulo $2m$).  One
can verify that $G$ is $3$--edge--connected. In
Figure \ref{3connesso}, you can see a graph $G$ of order $2n=240$.

We show that $\chi'_{[n-1]}(G)>4$. Suppose $\chi'_{[n-1]}(G)=4$ and
denote by $\{M_1, M_2, M_3, M_4\}$ a $[n-1]$-covering of $G$ of
size $4$. Since $2m\geq 10$, there is (at least) one copy of $H$, say
$H_1$, whose vertices are all covered in $M_j$, for every $1\leq
j\leq 4$. Furthermore, every $M_j$ contains a copy of exactly one
of the edges $e_i$'s, otherwise $H_1$ has uncovered vertices in $M_j$.
Therefore, we obtain an excessive factorization of size $4$ for each copy of the Petersen graph in $H_1$, a contradiction by $\chi'_{[n]}(P)=5$ again. 
\qed

\subsection{Oddness $2$ and $4$}

The oddness of a cubic graph $G$ is the minimum number of odd circuits
in a $2$-factor of $G$. Obviously, the oddness of a cubic graph is
an even number and it is $0$ if and only if the graph is
$3$--edge--colorable. The next results hold for cubic graph of
oddness $2$ and $4$.

\begin{proposition}\label{pro:oddness2}
Let $G$ be a cubic graph of order $2n$ and oddness $2$, then
$\chi'_{[n-1]}(G)=4$.
\end{proposition}

\textit{Proof.} Let $F$ be a $2$--factor of $G$ having exactly $2$
odd circuits. Let $M$ be the complementary perfect matching of $F$
in $G$. Let $N$ be a $[n-1]$-matching of $F$. Since $M$ and $N$ are
disjoint the assertion follows from Lemma \ref{MN}. $\qed$

\smallskip

As a consequence of Proposition \ref{pro:oddness2}, the excessive
$[n-1]$--index of every permutation snark of order $2n$ is $4$. We
recall that a permutation snark is a snark containing a $2$--factor
of exactly $2$ odd circuits having no chords. The Petersen graph is a
permutation snark.

\bigskip

In order to prove an equivalent result for cubic graphs of oddness $4$, we need the following lemma proved in \cite{KaiKraNor}. 
We refer to Schrijver's monography \cite{Sch} for the definition of fractional perfect matching and related topics.

\begin{lemma}\label{tool}
If $w$ is a fractional perfect matching in a cubic graph $G$ and $c \in \mathbb{R}^E$, then $G$ has a perfect matching $M$ such that $$ c \cdot \chi^M \geq c \cdot w $$
where $\cdot$ denotes the scalar product.
\end{lemma}

\begin{proposition}\label{pro:oddness4}
Let $G$ be a cyclically $4$-connected cubic graph of order $2n$ and oddness $4$. Then, $\chi'_{[n-1]}(G)=4$.
\end{proposition}
\textit{Proof.}
Let $M_1$ be a perfect matching of $G$ whose complementary $2$-factor $F$ has exactly $4$ odd circuits.
Since $G$ is cyclically $4$-connected, the function $\omega:E(G) \to \mathbb{R}$, defined by

$$ \omega(e)=\left\{ \begin{array}{ll}
                      \frac{1}{5} & \mbox{if } e \in M_1 \\
                      & \\
                      \frac 25 & \mbox{if } e \notin M_1 \\
                     \end{array}\right.
$$
for each $e \in E(G)$, is a fractional perfect matching of $G$ (see for instance \cite{KaiKraNor}).\\
Select a pair of incident edges in each odd circuit of $F$, say $L$ the set of these eight edges of $G$.
By Lemma \ref{tool}, there exists a perfect matching $M_2$ such that 
$$\chi^{L} \cdot \chi^{M_2} \geq \chi^{L} \cdot w.$$
The left hand side of the previous inequality is exactly the number of edges in $M_2 \cap L$ and the right one is equal to $8 \cdot \frac 25 = \frac{16}{5}$. Whence, $|M_2 \cap L|=4$.\\
The subgraph $H=G \setminus \{M_1 \cup (M_2 \cap L) \}$ consists of $4$ paths of even length and, in case, even cycles. Hence,
the edges of $H$ can be covered by two $[n-2]$-matchings, say $N_1$ and $N_2$.
The set $\{M_1, M_2, N_1, N_2\}$ satisfies the hypothesis of Lemma \ref{lemmacovering} and the assertion follows. $\qed$ \\

\subsection{$3^*$-connected graphs}

The class of $3^*$-connected cubic graphs is first considered in
\cite{AlbAldHolShe}. A $3$-graph is said to be $3^*$-connected if
there exists a pair of vertices $a,b \in V(G)$ such that $a,b$ are
the endvertices of three openly disjoint paths $Q_1,Q_2, Q_3$ such
that $V(G)=\bigcup_{i=1}^3 V(Q_i)$.

It is natural to give an equivalent definition of the class of
$3^*$-connected cubic graphs in our context. A $3^*$-connected cubic
graph is a $3$-graph having a $[n-1]$--matching $M$
such that $G \smallsetminus M$ is connected.\\
In the following proposition, we determine the excessive
$[n-1]$-index of a $3^*$-connected cubic graph.

\begin{proposition}\label{pro:3connected}
Let $G$ be a $3^*$-connected cubic graph of order $2n$, with
$n\geq4$. Then $\chi_{[n-1]}(G)=4$.
\end{proposition}

\textit{Proof.} By the definition of $3^*$-connected cubic graph,
there exist two vertices $a,b$ and three openly disjoint paths
$Q_1=(a,u_1,\ldots,u_r,b)$, $Q_2=(a,v_1,\ldots,v_s,b)$ and
$Q_3=(a,w_1,\ldots,w_t,b)$ spanning the graph $G$. Without loss of
generality, we can assume $r$ even and hence $s+t$ even. Let $N$ be
the complementary $[n-1]$-matching of $Q_1 \cup Q_2 \cup Q_3$.

Denote by $C$ the $(s+t+2)$--cycle $C=Q_2\cup Q_3$ and by $Q$ the
subpath $Q=(u_1, u_2,\ldots, u_r)$ of $Q_1$.

As $C$ has even length and  $Q$ has odd length (i.e. $Q$ has an odd
number of edges), we can color alternately the edges of $C$ and
those of $Q$ obtaining a perfect matching $M$ of $G$.

Since $M$ and $N$ are disjoint the assertion follows by Lemma
\ref{MN}. $\qed$

\subsection{Circumference}

The circumference of a graph $G$ is the length of any longest circuit
of $G$. In the next proposition we give a further support to the claim that
large classes of cubic graphs have excessive $[n-1]$-index equals to
$4$.

\begin{proposition}\label{pro:circumference}
Let $G$ be a cubic graph of order $2n\geq 8$ and circumference at least $2n-2$. Then $\chi'_{[n-1]}(G)=4$.
\end{proposition}

\textit{Proof.} If $G$ has circumference $2n$, then $G$ is
$3$--edge--colorable and the assertion follows from Proposition
\ref{class1}; if $G$ has circumference $2n-1$, then $G$ is
$3^*$--connected and the assertion follows from Proposition
\ref{pro:3connected}.

We consider $G$ with circumference $2n-2$. Denote by $C$ a circuit of
$G$ of length $2n-2$ and by $u$, $v$ the vertices of $G$ not
belonging to $C$. We color alternately the edges of $C$ and obtain
two $[n-1]$--matchings of $G$, say $M_1$ and $N$. We distinguish two
cases according that $u$, $v$ are adjacent or not.

If $u$, $v$ are adjacent vertices in $G$, then $M=M_1\cup\{[u, v]\}$
and $N$ satisfy Lemma \ref{MN} and the assertion follows.

Consider  $u$, $v$ non--adjacent. The set $L$ of chords of $C$ is a
$[n-4]$-matching of $G$ (the vertices of $C$ which are adjacent to
$u$ and $v$ are uncovered in $L$). Denote by $u_i$ (respectively, by
$v_i$) the vertices of $C$ adjacent to $u$ (respectively, to $v$),
with $i=1, 2, 3$.

The subgraph $M_1\cup L$ contains exactly three paths of odd length
whose endvertices are the uncovered vertices of $L$.  At least one
of these three paths has one endvertex adjacent to $u$, say $u_1$,
and the other adjacent to $v$, say $v_1$. The subgraph $H=M_1\cup L
\cup \{[u, u_i], [v, v_i] : i=1,2\}$ is a $2$--edge--colorable.
Since exactly one connected component of $H$ is a path of odd
length, a color class of $H$ is a perfect matching of $G$ and the
other one is a $[n-1]$--matching of $G$. Again, the assertion
follows from Lemma \ref{MN}.$\qed$

\section{Final Remarks}
The main aim of this paper has been the study of the excessive $[n-1]$-index of cubic graphs. Among other results,
we have given some evidence that each cyclically $4$-connected cubic graph of order $2n$ can be covered with four $[n-1]$-matchings. 
Nevertheless, analogously to perfect matchings case, we leave completely open the weaker problem of the existence of a constant $k$ such that $\chi'_{[n-1]}(G)$ is at most $k$ for every $3$-graph $G$.

\end{document}